\documentclass[11pt]{article}
\usepackage{amssymb,amsthm,amsmath,graphicx}

\newtheorem{theorem}{Theorem}[section]
\newtheorem{remark}{Remark}[section]
\newtheorem{proposition}{Proposition}[section]
\newtheorem{lemma}{Lemma}[section]
\numberwithin{equation}{section}

\newcommand{\h}{\hspace*{.24in}}
\def\geqslant {\ge}
\def\leqslant {\le}
\def\bqq{\begin{eqnarray*}}
\def\eqq{\end{eqnarray*}}
\def\bq{\begin{eqnarray}}
\def\eq{\end{eqnarray}}
\def\nn{\nonumber}

\title{\bf Ice formation in the Arctic during summer: false-bottoms}
\author{Phan Thanh Nam$^{a,b}$, Pham Ngoc Dinh Alain$^a$\thanks{Author correspondence:
alain.pham@univ-orleans.fr}\\
{Pham Hoang Quan$^{b}$ and Dang Duc Trong$^b$}
\\\\\small\it$^a$Mathematics Department, MAPMO UMR 6628, BP 67-59, 45067 Orleans cedex, France
\\\small\it$^b$Department of Mathematics , HoChiMinh City National University, Viet Nam}
\date{{}}

\begin{document}
\maketitle
\begin{abstract} The only source of ice formation in the Arctic during summer is a layer of ice called false-bottoms
between an under-ice melt pond and the underlying ocean. Of interest is to give a mathematical model in order
to determine the simultaneous growth and ablation of false-bottoms, which is governed by both of heat fluxes and salt fluxes.
In one dimension, this problem may be considered mathematically as a two-phase Stefan problem with two free boundaries.
Our main result is to prove the existence and uniqueness of the solution from the initial condition.
\end{abstract}

\section{Introduction}\label{sec1}

\h There are some different reservoirs of fresh water in the Arctic during summer (see, e.g., {\it Eicken et al.} \cite{EKKP}). First, melt water collects in surface melt pond (melting under the sun) which is the most important reservoir. Second, this melt water can percolate into the ice matrix to form an {\it under-ice melt pond} (see \cite{H} for more detail). At the interface between this fresh water and the underlying salt water, double-diffusive convection of heat and salt leads to the formation of a layer of ice called {\it false-bottoms} (see Figure 1 below). Very early, {\it Nansen} \cite{N} in 1897 noted that this is the only source of forming new ice in the Arctic during the summer. This phenomenon has been considered for a long time by many authors (see, e.g. \cite{E, EKKP,H, MK,MMN, NMWMSE,UB}). However, it has been considered in geophysical view-point based on practical experiments rather than rigorously mathematical formulations.

One of the most interesting ones is the simultaneous growth and
ablation of false-bottoms, which is governed by both of heat fluxes
and salt fluxes. The ablation of the sea-ice interface is caused by
dissolution rather than by melting. Note that salt water has the
double properties:  it does not freeze even for temperature less
than $0^0$C, and it dissolves ice when it is in contact with ice.
{\it McPhee et al.} \cite{MMN} emphasized that properly describing
heat and salt flux at the ice-ocean interface is essential for
understanding and modeling the false-bottoms, and in particular
without the double diffusion at this interface false bottoms would
be so short-lived. The growth of the upper interface between a false
bottom and a under-ice melt pond is governed by the purely
thermodynamic condition at the interface.
\newpage
\centerline{\includegraphics[width=7cm,height=7cm,
keepaspectratio=true]{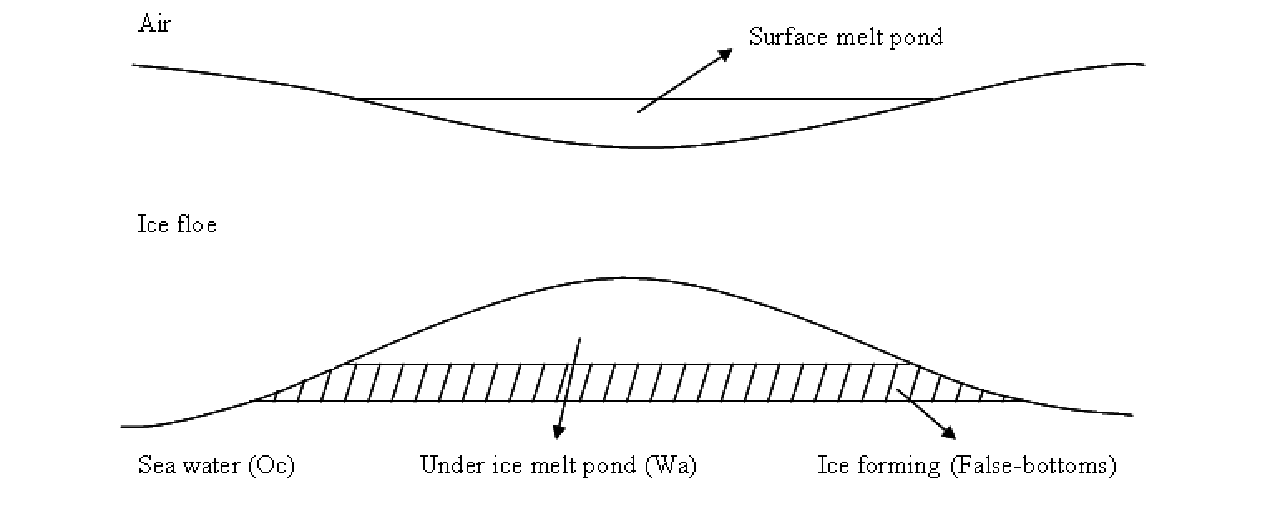}}
\begin{center}{Figure 1: Ice formation in the Arctic during summer.}\end{center}

Recently, in 2003, {\it Notz et al.} \cite{NMWMSE} gave a model simulating successfully the simultaneous growth and ablation of false-bottoms.  They formulated mathematically the problem by a system of partial differential equations and solved them numerically by using a numerical routine in Mathematica. Although their numerical result fits quite well to early experimental data from {\it Martin and Kauffman} \cite{MK}, a rigorous proof of the existence and uniqueness of the solution for the system of equations is still unavailable. Our aim in this paper is to give  a such a mathematical proof. More precisely, we shall represent the problem explicitly by a system of partial differential equations associated with free boundary conditions similar to \cite{NMWMSE}, and then show that the system has a unique solution from given initial conditions.

Now we consider a one-dimensional model describing the simultaneous growth and ablation of the ice of false-bottoms. Here we have three environments: the ocean (Oc), the ice of false-bottoms (Fb) and the fresh water (Wa).  Denote by $T(x,t)$, $S(x,t)$ the temperature and the salinity, and denote by $h_0(t),h_u(t)$ the free boundaries at the interfaces ice-ocean (Fb-Oc) and ice-water (Fb-Wa), respectively (see Figure 2 below).

\centerline{\includegraphics[width=7cm,height=7cm,
keepaspectratio=true]{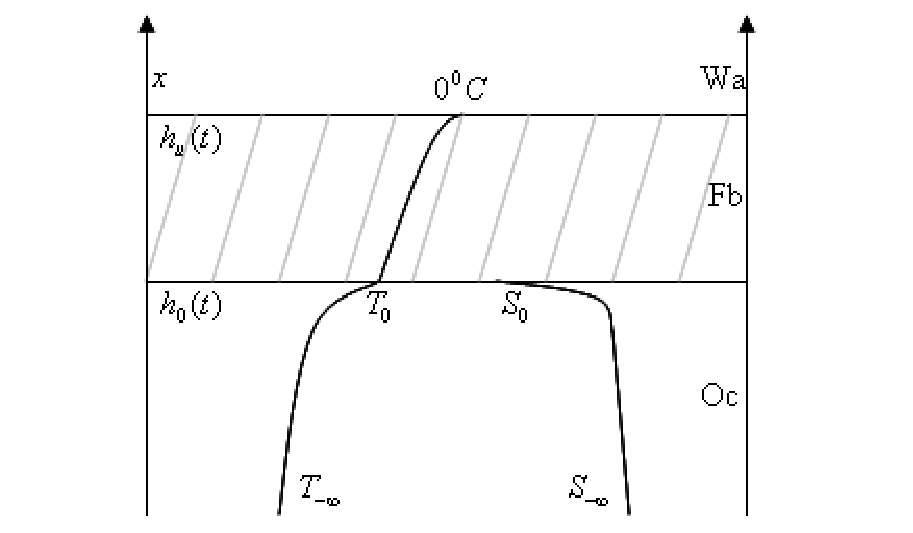}}
\begin{center}\text{Figure 2: One-dimensional model.}\end{center}

At the interface Fb-Oc, we apply the first principle of thermodynamics, (i.e. variation of energy = variation of heat flux through the interface)
\[
\Delta U = \Delta \Phi .
\]
The net amount of heat transferred through the interface Fb-Oc  in a section s is equal to
$$\Delta U=-dh_0\rho_IL_fs,$$
where $\rho_I$ is the density of the ice and $L_f$ is the latent heat of fusion. On the other hand, the difference of the heat fluxes through a section s in the ice and the ocean during a time $dt$ is
$$
\Delta \Phi=\left( {-\lambda _I T_x(h_0(t)+,t)  + \lambda _O T_x(h_0(t)-,t)} \right)sdt,
$$
where $\lambda _I$, $\lambda _O$ are thermal conductivities of the ice and the ocean. Here the notations $h_0(t)+$ and $h_0(t)-$ stand for the right limit and the left limit at $x=h_0(t)$. Thus the law of conservation of energy mentioned above, i.e. $\Delta U = \Delta \Phi$, leads to the Stefan condition for the heat balance at the interface
\bq
h_0 '(t) = \widetilde\lambda _I T_x(h_0(t)+,t)  - \widetilde\lambda _O T_x(h_0(t)-,t), \label{Fb-Oc1}
\eq
where
$$
\widetilde\lambda _I  = \frac{{\lambda _I }}
{{\rho _I L_f }}>0,~\widetilde\lambda _O  = \frac{{\lambda _O }}
{{\rho _I L_f }}>0.
$$
For simplicity, we can neglect the salt of the ice of false-bottoms. The water near the interface Fb-Oc is a mixture of melt water, which melts from the ice of false-bottoms, and sea water. This water freshens at the rate $S_0(t)h_0'(t)$, while salt diffuses into this water at the rate $- DS_x(h_0(t)-,t)$,
where $S_0(t)=S(h_0(t)-,t)$ is the salinity of the ocean at the interface and $D>0$ is the molecular diffusivity of salt in sea water. The balance of salt at this interface leads to the conservation condition
\bq
S_0(t)h_0 '(t) =  - DS_x(h_0(t)-,t).  \label{Fb-Oc2}
\eq
Moreover, the interface temperature $T_0 = T(h_0(t),t)$ and salinity $S_0$ are connected via freezing-point relationship
\bq
T_0+n_0T_0'=-m_0S_0,  \label{T0S0}
\eq
where $m_0,n_0$ are non-zero constants. The establishment of equations (\ref{Fb-Oc1}) and (\ref{Fb-Oc2}) bases on {\it Martin and Kauffman} \cite{MK} and {\it Notz et al.} \cite{NMWMSE}.  In (\ref{T0S0}), if $n_0=0$ and $m_0=0.054^0\text{C psu}^{-1}$ then we obtain the approximation $T_0\approx -m_0S_0$, which was used in \cite{MK, NMWMSE}. However, as we shall show later, it is reasonable to assume that the salinity is dependent on the gradient of the temperature at the freezing point.

At the interface Fb-Wa, we use the simplified scheme presented by {\it Grenfell and Maykut} \cite{GM}, in which the temperature of the water in the under-ice melt pond is kept at $0^0$C. In particular, it leads to the boundary condition at the interface
\bq
T(h_u (t)-,t) = 0.\label{Fb-Wa1}
\eq
Furthermore, due to the neglect of heat fluxes from the fresh water above the false-bottoms, a thermodynamic condition similar to (\ref{Fb-Oc1}) reduces to the Stefan condition at the upper surface
\bq
h_u '(t) = \widetilde\lambda _I T_x(h_u(t)-,t).\label{Fb-Wa}
\eq

Finally, we have the following diffusion equations for salt and heat in the ocean and in the ice
\bq
T_t  &=& D_I T_{xx},\h h_u(t)>x>h_0(t),\label{Fb}\hfill\\
T_t  &=& D_O T_{xx} ,\h h_0(t)>x>-\infty,\label{OcT}\hfill\\
S_t&=& DS_{xx},\h h_0(t)>x>-\infty,\label{OcS}
\eq
where $D_I>0$ and $D_O>0$ are the thermal diffusivity. Here all mentioned constants are given.

Assume that the initial conditions $h_0(0)=h^0_0$, $h_u(0)=h^0_u$, $T(x,0)=T^0(x)$ in $h_u(0)>x>-\infty$ and $S(x,0)=S^0(x)$ in $h_0(0)>x>-\infty$ are given. The problem is of finding from the initial conditions $(h_0^0,h^0_u,T^0,S^0)$ a solution $(h_0,h_u,T,S)$ of the system (\ref{Fb-Oc1})-(\ref{OcS}), where $T=T(x,t)$ in $t>0$, $h_u(t)>x>-\infty$ and $S=S(x,t)$ in $t>0$, $h_0(t)>x>-\infty$. This is a one-dimensional two-phase Stefan problem with two free boundaries.

The main result of the present paper is to prove the existence and uniqueness a local solution for the problem. Rigorously, we suppose that the initial conditions $(h_0^0,h_u^0,T^0,S^0)$ satisfy the following three assumptions

(H1) $h_0^0<h_u^0$;

(H2) $(T^0)_x$ is continuous and bounded in $x\in (-\infty, h_0^0]$ and $x\in [h_0^0, h_u^0]$; $T^0$ is continuous at $x=h_0^0$ and $T^0(h_u^0-)=0$;

(H3) $S^0$ is continuous and bounded in $x\in (-\infty, h_0^0]$.

\begin{remark} The condition (H1) means the ice layer of the false-bottom has already existed. It is of course a much more interesting problem to consider the behavior of the false-bottom at the starting time, which corresponds to condition $h_0^0=h_u^0$, but it is extremely difficult. We mention that the difficulty does not only come from mathematical computation but also be recognized by physical experiments. For example, {\it Notz et al.} \cite{NMWMSE} simulated the experiment of {\it Martin and Kauffman} \cite{MK}, in which they put fresh water at $0^0$C on top of salt water at its freezing point in order to simulate the evolution of a false-bottoms, but this model does not take salt transport through the false-bottom into account. Therefore they then started with a 5-cm layer of ice.
\end{remark}

\begin{remark} It is implicit from (\ref{Fb-Oc1}) and (\ref{Fb-Wa}) that $T_x(.,t)$ is continuous at $x=h_0(t)-$, $h_0(t)+$ and $h_u(t)-$. Therefore, it is reasonable to assume that $(T^0)_x$ is continuous also at $x=h_0^0-$, $h_0^0+$ and $h_u^0-$ in (H2). Similarly, it is natural to assume that $T^0$ is continuous at $x=h_0^0$ and $T^0(h_u^0-)=0$ due to (\ref{T0S0}) and (\ref{Fb-Wa1}).
\end{remark}

Let $\sigma>0$. We say that $(h_0,h_u,T,S)$ is a solution of the system (\ref{Fb-Oc1})-(\ref{OcS}) in $(0,\sigma)$ corresponding to the initial conditions $(h_0^0, h_u^0,T^0,S^0)$ if the following four conditions holds

(C1) $h_0(t)$ and $h_u(t)$ are continuously differentiable in $t\in [0,\sigma)$; $h_0(0)=h_0^0<h_u^0=h_u(0)$ and $h_0(t)<h_u(t)$ in $t\in (0,\sigma)$;

(C2) $T_t$, $T_{xx}$ is continuous in $t\in (0,\sigma)$, $x\in (-\infty,h_0(t))\cup (h_0(t),h_u(t))$; $T_{x}(.,t)$ is continuous at $x=h_0(t)-$, $x=h_0(t)+$ and $x=h_u(t)-$; $T(.,t)$ is continuous at $x=h_0(t)$; $T(x,.)$ is continuous at $t=0+$ and $T(x,0+)=T^0(x)$ in $x\in (-\infty,h_u^0)$;

(C3) $S_t$, $S_{xx}$ is continuous in $t\in (0,\sigma)$, $x\in (-\infty,h_0(t))$; $S_x(.,t)$ is continuous at $x=h_0(t)-$; $S(x,.)$ is continuous at $t=0$ and $S(x,0+)=S^0(x)$ in $x\in (-\infty,h_0^0)$;

(C4) Equations (\ref{Fb-Oc1})-(\ref{OcS}) hold in $t\in (0,\sigma)$.

In general, the conditions (C1)-(C4) assure a classical solution for the system of differential equations (\ref{Fb-Oc1})-(\ref{OcS}). Our main result is as follows.

\begin{theorem}\label{existenceanduniqueness} Assume that $(h_0^0, h_u^0,T^0,S^0)$ satisfy (H1)-(H3). Then there is a unique solution $(h_0,h_u,T,S)$ of the system (\ref{Fb-Oc1})-(\ref{OcS}) in $(0,\sigma)$ corresponding to the initial conditions $(h_0^0, h_u^0,T^0,S^0)$ for some $\sigma>0$. Moreover, this solution can be extended uniquely whenever the condition $h_u(\sigma)>h_0(\sigma)$ still holds.
\end{theorem}

Our proof follows the approach of A.Friedman \cite{F} (Chapter 8), which dealt with a classical one-phase Stefan problem with one free boundary. More precisely, we first reduce the problem to solving a system of nonlinear Volterra integral equations of the second kind and then solve this system by contraction principle. The remainder of the paper is divided into there sections. Section 2 is devoted to some preliminary results on Green functions and Volterra integral equations. In Section 3 we shall reformulate the problem to a system of nonlinear integral equations. In Section 4 we apply the contraction principle to prove the existence and uniqueness of a local solution for this system.

\section{Preliminaries}\label{pre4}

\h Let $a(t),b(t)$ be continuously differentiable functions and $b(t)>a(t)$ for all $t\ge 0$. Let $\kappa>0$ be a constant and let $u(x,t)$ be a solution of the diffusion equation
\bq
\frac{{\partial u}}
{{\partial t}} = \kappa \frac{{\partial ^2 u}}
{{\partial x^2 }},~~t>0,~b(t)> x>a(t).\label{heat}
\eq
We introduce the Green's function of equation (\ref{heat}),
\[
G(x,t;\xi ,\tau ) = \frac{H(t-\tau)}
{{2\sqrt {\pi \kappa (t - \tau )} }}\exp \left( { - \frac{{(x - \xi )^2 }}
{{4\kappa(t - \tau )}}} \right),
\]
where $H$ is the Heaviside function,
$$
H(t) = \left\{ \begin{gathered}
  1, ~\text{if}~t>0,\hfill \\
  0,~\text{if}~t<0. \hfill \\
\end{gathered}\right.
$$
The following lemma is useful to reformulate the differential equation (\ref{heat}) to a Volterra integral equations of the second kind.

\begin{lemma}\label{mainLemma} If $u$ is a solution of (\ref{heat}) then for $t>0$ and $a(t)<x<b(t)$ we have
\bqq
  u(x,t) &=& \int\limits_0^t {G(x,t;b(\tau ),\tau )\left[ {\kappa u_\xi  (b(\tau ) - ,\tau ) + u(b(\tau ),\tau )b'(\tau )} \right]d\tau }  \hfill \\
   &~&- \int\limits_0^t {G(x,t;a(\tau ),\tau )\left[ {\kappa u_\xi  (a(\tau ) + ,\tau ) + u(a(\tau ),\tau )a'(\tau )} \right]d\tau }  \hfill \\
   &~&- \kappa \int\limits_0^t {G_\xi  (x,t;b(\tau ),\tau )u(b(\tau ),\tau )d\tau }  + \kappa \int\limits_0^t {G_\xi  (x,t;a(\tau ),\tau )u(a(\tau ),\tau )d\tau }  \hfill \\
  &~&+\int\limits_{a(0)}^{b(0)} {G(x,t;\xi ,0)u(\xi ,0)d\xi } .
\eqq
\end{lemma}
\begin{proof} Note that
\bq
G_\tau   + \kappa G_{\xi \xi }  = 0 ~\text{for all}~ \tau<t,~~\text{and} ~G(x,t;\xi ,t - ) = \delta (x - \xi ),
\eq
where $\delta=H'$ is Dirac delta function.

Integrating the Green's identity, here $u=u(\xi,\tau)$,
\[
\kappa\frac{\partial }
{{\partial \xi }}\left( {G\frac{{\partial u}}
{{\partial \xi }} - u\frac{{\partial G}}
{{\partial \xi }}} \right) - \frac{\partial }
{{\partial \tau }}\left( {uG} \right) = 0
\]
over the domain $a(\tau ) < \xi  < b(\tau ),0 < \tau  < t$, we will obtain the desired result because
\bqq
&~&\int\limits_0^t {\int\limits_{a(\tau )}^{b(\tau )} {\frac{\partial }
{{\partial \xi }}\left( {G\frac{{\partial u}}
{{\partial \xi }} - u\frac{{\partial K}}
{{\partial \xi }}} \right)d\xi d\tau } }  = \int\limits_0^t {\left[ {G\frac{{\partial u}}
{{\partial \xi }} - u\frac{{\partial G}}
{{\partial \xi }}} \right]_{\xi  = a(\tau )}^{\xi  = b(\tau )} d\tau }  \hfill \\
   &=& \int\limits_0^t {G(x,t;b(\tau ),\tau )u_\xi  (b(\tau ) - ,\tau )d\tau }  - \int\limits_0^t {G(x,t;a(\tau ),\tau )u_\xi  (a(\tau ) + ,\tau )d\tau }  \hfill \\
   &~&- \int\limits_0^t {G_\xi  (x,t;b(\tau ),\tau )u(b(\tau ),\tau )d\tau }  + \int\limits_0^t {G_\xi  (x,t;a(\tau ),\tau )u(a(\tau ),\tau )d\tau } ,
\eqq
and
\bqq
&~&\int\limits_0^t {\int\limits_{a(\tau )}^{b(\tau )} {\frac{\partial }
{{\partial \tau }}\left( {uG} \right)d\xi d\tau } }\hfill \\
   &=& \int\limits_0^t {\left\{ {\frac{\partial }
{{\partial \tau }}\left( {\int\limits_{a(\tau )}^{b(\tau )} {uGd\xi } } \right) - [uG]_{\xi  = b(\tau )} b'(\tau ) + [uG]_{\xi  = a(\tau )} a'(\tau )} \right\}d\tau }  \hfill \\
   &=& \left[ {\int\limits_{a(\tau )}^{b(\tau )} {uGd\xi } } \right]_{\tau  = 0}^{\tau  = t - }  - \int\limits_0^t {[uG]_{\xi  = b(\tau )} b'(\tau )d\tau }  + \int\limits_0^t {[uG]_{\xi  = a(\tau )} a'(\tau )d\tau }  \hfill \\
   &=& u(x,t) - \int\limits_{a(0)}^{b(0)} {G(x,t;\xi ,0)u(\xi ,0)d\xi }  - \int\limits_0^t {G(x,t;b(\tau ),\tau )u(b(\tau ),\tau )b'(\tau )d\tau }  \hfill \\
 &~&  + \int\limits_0^t {G(x,t;a(\tau ),\tau )u(a(\tau ),\tau )a'(\tau )d\tau } .
\eqq
Here we have made use the following equality, with $v=uG$,
$$
\int\limits_{a(\tau )}^{b(\tau )} {\frac{\partial }
{{\partial \tau }}\left( {v(\xi ,\tau )} \right)d\xi }  = \frac{\partial }
{{\partial \tau }}\left( {\int\limits_{a(\tau )}^{b(\tau )} {v(\xi ,\tau )d\xi } } \right) - v(b(\tau ),\tau )b'(\tau ) + v(a(\tau ),\tau )a'(\tau ).
$$
In fact, in the case $a(t)\equiv 0$ we have
\bqq
&~& \int\limits_0^{b(\tau )} {\frac{\partial }
{{\partial \tau }}\left( {v(\xi ,\tau )} \right)d\xi }  = \int\limits_0^\infty  {H(b(\tau ) - \xi )\frac{\partial }
{{\partial \tau }}\left( {v(\xi ,\tau )} \right)d\xi }  \hfill \\
   &=& \int\limits_0^\infty  {\frac{\partial }
{{\partial \tau }}\left[ {H(b(\tau ) - \xi )v(\xi ,\tau )} \right]d\xi }  - \int\limits_0^\infty  {\frac{\partial }
{{\partial \tau }}\left[ {H(b(\tau ) - \xi )} \right]v(\xi ,\tau )d\xi }  \hfill \\
   &=& \frac{\partial }
{{\partial \tau }}\left( {\int\limits_0^\infty  {H(b(\tau ) - \xi )v(\xi ,\tau )d\xi } } \right) - \int\limits_0^\infty  {\delta (b(\tau ) - \xi )b'(\tau )v(\xi ,\tau )d\xi }  \hfill \\
   &=& \frac{\partial }
{{\partial \tau }}\left( {\int\limits_0^{b(\tau )} {v(\xi ,\tau )d\xi } } \right) - v(b(\tau ),\tau )b'(\tau ).
\eqq
If $a(t)$ is not constant, we can write
\[
\int\limits_{a(\tau )}^{b(\tau )} {\frac{\partial }
{{\partial \tau }}\left( {v(\xi ,\tau )} \right)d\xi }  = \int\limits_0^{b(\tau )} {\frac{\partial }
{{\partial \tau }}\left( {v(\xi ,\tau )} \right)d\xi }  - \int\limits_0^{a(\tau )} {\frac{\partial }
{{\partial \tau }}\left( {v(\xi ,\tau )} \right)d\xi } ,
\]
and apply the above result. This completes the proof.
\end{proof}
\begin{remark} The result in Lemma \ref{mainLemma} still holds for $b(t)\equiv +\infty$ or $a(t)\equiv -\infty$. For example, if $a(t)\equiv -\infty$ then the formula in Lemma \ref{mainLemma} reduces to
\bqq
  u(x,t) &=& \int\limits_0^t {G(x,t;b(\tau ),\tau )\left[ {\kappa u_\xi  (b(\tau ) - ,\tau ) + u(b(\tau ),\tau )b'(\tau )} \right]d\tau }  \hfill \\
   &~&- \kappa \int\limits_0^t {G_\xi  (x,t;b(\tau ),\tau )u(b(\tau ),\tau )d\tau }  + \int\limits_{ - \infty }^{b(0)} {G(x,t;\xi ,0)u(\xi ,0)d\xi } .
\eqq
\end{remark}

We shall need also an useful lemma giving the jump relation at the boundary (see Friedman \cite{F}, page 217, Lemma 1).

\begin{lemma}\label{jumpLemma} Let $p(t)$ be continuous and let $s(t)>0$ satisfy the Lipschitz condition, $0\le t\le \sigma$. Then, for $0< t\le \sigma$,
\[
\mathop {\lim }\limits_{x \to s(t)^ -  } \kappa\int\limits_0^t {p(\tau)G_x  (x,t;s(\tau ),\tau )d\tau }  = \frac{1}
{2}p(t) +\kappa\int\limits_0^t {p(\tau)G_x (s(t),t;s(\tau ),\tau )d\tau } .
\]
\end{lemma}

\begin{remark} In applications of the above lemma, sometimes we need to note that $G_{\xi}=-G_x$. Moreover, for the right limit we have
 \[
\mathop {\lim }\limits_{x \to s(t)^ +  } \kappa\int\limits_0^t {p(\tau)G_x  (x,t;s(\tau ),\tau )d\tau }  = -\frac{1}
{2}p(t) +\kappa\int\limits_0^t {p(\tau)G_x (s(t),t;s(\tau ),\tau )d\tau } .
\]
\end{remark}

Finally, we state a simple version of the uniqueness for a system of linear Volterra integral equations of the second kind.

\begin{lemma}\label{linearVolterra} Let $n \in \mathbb{N}$, $\sigma>0$, $p>1$, $q>1$, $1/p+1/q=1$. Assume that $\tau\mapsto W_j (t,\tau )$ is measurable in $(0,t)$ for all $t\in (0,\sigma]$ and
\[
\int\limits_0^t {\left| {W_j (t,\tau )} \right|^p d\tau }  \leqslant const.,~\forall t \in (0,\sigma ],j = \overline {1,n} .
\]
Then the system
\[
\Psi _i (t) = \sum\limits_{j = 1}^n {\left( {\int\limits_0^t {W_j (t,\tau )\Psi _j (\tau )d\tau } } \right)} ,\forall t \in (0,\sigma ],j = \overline {1,n},
\]
has a unique solution $\{\Psi _j\}_{j=1}^{n}=0$ in $L^q(0,\sigma)$.
\end{lemma}
\begin{proof} Using Holder inequality one has
\bqq
  \left| {\Psi _i (t)} \right| &\leqslant& \sum\limits_{j = 1}^n {\left| {\int\limits_0^t {W_j (t,\tau )\Psi _j (\tau )d\tau } } \right|}  \hfill\\
&\leqslant& \sum\limits_{j = 1}^n {\left( {\int\limits_0^t {\left| {W_j (t,\tau )} \right|^p d\tau } } \right)^{1/p} \left( {\int\limits_0^t {\left| {\Psi _j (\tau )} \right|^q d\tau } } \right)} ^{1/q}  \hfill \\
   &\leqslant& const.\sum\limits_{j = 1}^n {\left( {\int\limits_0^t {\left| {\Psi _j (\tau )} \right|^q d\tau } } \right)} ^{1/q}.
\eqq
Therefore,
\[
\sum\limits_{j = 1}^n {\left| {\Psi _j (\tau )} \right|^q }  \leqslant const.\sum\limits_{j = 1}^n {\left( {\int\limits_0^t {\left| {\Psi _j (\tau )} \right|^q d\tau } } \right)}  = const.\int\limits_0^t {\left( {\sum\limits_{j = 1}^n {\left| {\Psi _j (\tau )} \right|^q } } \right)} d\tau ,
\]
and it follows from Gronwall's Lemma that $\sum\limits_{j = 1}^n {\left| {\Psi _j (\tau )} \right|^q }=0$. Thus $\{\Psi _j\}_{j=1}^{n}=0$.
\end{proof}
\begin{remark} If , for example,
\[
\left| {W_j (t,\tau )} \right| \leqslant \frac{{const.}}
{{\sqrt {t - \tau } }},~~j = \overline {1,n} .
\]
then we can choose any $p\in (1,2)$ in order to apply Lemma \ref{linearVolterra}.
\end{remark}

\section{Reduction to integral equations}

\h Denote by $G_1, G_2 ,G_3$ the Green function $G$ in Lemma \ref{mainLemma} corresponding to $\kappa=D,D_O ,D_I$, respectively. We shall reformulate the problem (\ref{Fb-Oc1})-(\ref{OcS}) to a system of integral equations of time-depending functions $(v_0,v_1,v_2,v_3)$ where
\bq
v_0 (t) = T_0 '(t),~v_1 (t) = T_x (h_0 (t) - ,t),~v_2 (t) = T_x (h_0 (t) + ,t),~v_3 (t) = T_x (h_u (t) - ,t).\label{v}
\eq
Assume that $(h_0,h_u,T,S)$ is a solution of the system (\ref{Fb-Oc1})-(\ref{OcS}) in $(0,\sigma)$ corresponding to the initial conditions $(h_0^0, h_u^0,T^0,S^0)$ for some $\sigma>0$. Applying Lemma \ref{mainLemma} to equation (\ref{OcS}) with $h_0(t)>x>-\infty$ and using condition (\ref{Fb-Oc2}), one has
\bq
S(x,t) =  - D\int\limits_0^t {G_{1\xi}  (x,t;h_0 (\tau ),\tau )S_0 (\tau )d\tau }  + \int\limits_{ - \infty }^{h_0^0} {G_1(x,t;\xi ,0)S^0(\xi)d\xi } .\label{ieS}
\eq
In particular, we see from (\ref{ieS}) that $S(x,t)$ is determined completely by $h_0$ and $S_0(t)$. Taking $x\to h_0(t)-$ in (\ref{ieS}) and using the jump relation in Lemma \ref{jumpLemma}, we get
\bq
S_0 (t) =  - 2D\int\limits_0^t {G_{1\xi}  (h_0 (t),t;h_0 (\tau ),\tau )S_0 (\tau )d\tau }  + 2\int\limits_{ - \infty }^{h_0^0} {G_1(h_0 (t),t;\xi ,0)S^0(\xi)d\xi } .\label{ieS0}
\eq
Note that for $v_0 (t) = T_0 '(t)$ the condition (\ref{T0S0}) can be rewritten as
\bq
T^0 (h_0^0 ) + \int\limits_0^t {v_0 (\tau )d\tau }  + n_0 v_0 (t) =  - m_0 S_0 (t).\label{v0S0}
\eq
Here we have used $T_0(0)=T(h_0(0),0)=T_0(h_0^0)$. We deduce from (\ref{ieS0}) and (\ref{v0S0}) the first integral equation
\bq
  v_0 (t) &=&  - \frac{1}
{{n_0 }}T^0 (h_0^0 ) - \frac{1}
{{n_0 }}\int\limits_0^t {v_0 (\tau )d\tau }  - \frac{{2m_0 }}
{{n_0 }}\int\limits_{ - \infty }^{h_0^0 } {G_1 (h_0 (t),t;\xi ,0)S^0 (\xi )d\xi }  \nn\hfill \\
   &~&- \frac{{2D}}
{{n_0 }}\int\limits_0^t {G_{1\xi } (h_0 (t),t;h_0 (\tau ),\tau )\left[ {T^0 (h_0^0 ) + \int\limits_0^\tau  {v_0 (s)ds}  + n_0 v_0 (\tau )} \right]d\tau } .\label{iev0}
\eq

We next apply Lemma \ref{mainLemma} to equation (\ref{OcT}) for $h_0(t)>x>-\infty$ to get
\bq
  T(x,t) &=& \int\limits_0^t {G_2 (x,t;h_0 (\tau ),\tau )\left[ {D_O v_1(\tau ) + T_0 (\tau )h_0 '(\tau )} \right]d\tau }\nn  \hfill \\
 &~& - D_O \int\limits_0^t {G_{2\xi } (x,t;h_0 (\tau ),\tau )T_0 (\tau )d\tau }  + \int\limits_{ - \infty }^{h_0^0} {G_2 (x,t;\xi ,0)T^0 (\xi )d\xi } ,  \label{ieTOc}
\eq
where $v_1(t)=T_x(h_0(t)-,t)$. We now differentiate both sides of (\ref{ieTOc}) with respect to $x$, then take $x\to h_0(t)-$. To go into the details, because $D_O{G_{2\xi x} }=-D_O{G_{2\xi \xi} }={G_{2\tau} }$, we have
\[
\begin{gathered}
~~   - D_O \int\limits_0^t {G_{2\xi x} (x,t;h_0 (\tau ),\tau )T_0 (\tau )d\tau }  =  - \int\limits_0^t {G_{2\tau } (x,t;h_0 (\tau ),\tau )T_0 (\tau )d\tau }  \hfill \\
   =  - \int\limits_0^t {\left[ {\frac{\partial }
{{\partial \tau }}\left( {G_2 (x,t;h_0 (\tau ),\tau )} \right) - G_{2\xi } (x,t;h_0 (\tau ),\tau )h_0 '(\tau )} \right]T_0 (\tau )d\tau }  \hfill \\
   = G_2 (x,t;h_0^0,0)T_0 (0) + \int\limits_0^t {G_2 (x,t;h_0 (\tau ),\tau )T_0 '(\tau )d\tau }  \hfill \\
 ~~  - \int\limits_0^t {G_{2x} (x,t;h_0 (\tau ),\tau )T_0 (\tau )h_0 '(\tau )d\tau } . \hfill \\
\end{gathered}
\]
Moreover,
\bqq
&~& \int\limits_{ - \infty }^{h_0^0} {G_{2x} (x,t;\xi ,0)T^0 (\xi )d\xi }  =  - \int\limits_{ - \infty }^{h_0^0} {G_{2\xi } (x,t;\xi ,0)T^0 (\xi )d\xi }  \hfill \\
   &=&  - G_2 (x,t;h_0^0,0)T^0 (h_0^0) + \int\limits_{ - \infty }^{h_0^0} {G_2 (x,t;\xi ,0)T_\xi ^0 (\xi )d\xi }.
\eqq
Thus it follows from (\ref{ieTOc}) that
\bq
  T_x (x,t) &=& D_O \int\limits_0^t {G_{2x} (x,t;h_0 (\tau ),\tau )v_1(\tau )d\tau }\nn\hfill\\
&~& + \int\limits_0^t {G_2 (x,t;h_0 (\tau ),\tau )v_0(\tau )d\tau } +\int\limits_{ - \infty }^{h_0^0} {G_2 (x,t;\xi ,0)T_\xi ^0 (\xi )d\xi }\label{TxOc}
\eq
for all $h_0(t)>x>-\infty$. Here we have used the compatible condition $T_0(0)=T(h_0(0),0)=T^0(h_0^0)$ and replaced $T_0'(\tau)$ by $v_0(\tau)$. Taking $x\to h_0(t)-$ in (\ref{TxOc}) and using Lemma \ref{mainLemma} for the first term of the right hand side, we have the second integral equation

\bq
  v_1 (t) &=& 2D_O \int\limits_0^t {G_{2x} (h_0 (t),t;h_0 (\tau ),\tau )v_1 (\tau )d\tau }\nn\hfill\\
&~& + 2\int\limits_0^t {G_2 (h_0 (t),t;h_0 (\tau ),\tau )v_0(\tau )d\tau } +  2\int\limits_{ - \infty }^{h_0^0} {G_2 (h_0 (t),t;\xi ,0)T_\xi ^0 (\xi )d\xi }. \label{iev1}
\eq

Now we consider the heat distribution in false-bottoms. Applying Lemma \ref{mainLemma} to equation (\ref{Fb}) for $h_u(t)>x>h_0(t)$ one has
\bq
  T(x,t) &=& D_I \int\limits_0^t {G_3 (x,t;h_u (\tau ),\tau )v_3(\tau )d\tau }\nn  \hfill \\
&~&  - \int\limits_0^t {G_3 (x,t;h_0 (\tau ),\tau )\left[ {D_Iv_2(\tau ) + T_0 (\tau )h_0 '(\tau )} \right]d\tau }  \nn\hfill \\
&~&+ D_I \int\limits_0^t {G_{3\xi } (x,t;h_0 (\tau ),\tau )T_0 (\tau )d\tau }  + \int\limits_{h_0^0}^{h_u^0} {G_3 (x,t;\xi ,0)T^0(\xi)d\xi },\label{ieTFb}
\eq
where $v_2(t) = T_x (h_0 (t)+ ,t)$ and $v_3(t)= T_x (h_u (t)- ,t) = \dfrac{1}{{\widetilde{\lambda _I }}}h_u '(t)$. Let us differentiate both sides of (\ref{ieTFb}) with respect to $x$. We have
\bqq
&~& D_I \int\limits_0^t {G_{3\xi x} (x,t;h_0 (\tau ),\tau )T_0 (\tau )d\tau }  = \int\limits_0^t {G_{3\tau } (x,t;h_0 (\tau ),\tau )T_0 (\tau )d\tau }  \hfill \\
   &=& \int\limits_0^t {\left[ {\frac{\partial }
{{\partial \tau }}\left( {G_3 (x,t;h_0 (\tau ),\tau )} \right) - G_{3\xi } (x,t;h_0 (\tau ),\tau )h_0 '(\tau )} \right]T_0 (\tau )d\tau }  \hfill \\
   &=&  - G_3 (x,t;h_0^0,0)T_0 (0) - \int\limits_0^t {G_3 (x,t;h_0 (\tau ),\tau )T_0 '(\tau )d\tau }  \hfill \\
&~&  + \int\limits_0^t {G_{3x} (x,t;h_0 (\tau ),\tau )T_0 (\tau )h_0 '(\tau )d\tau }
\eqq
and
\bqq
  &~&\int\limits_{h_0^0}^{h_u^0} {G_{3x} (x,t;\xi ,0)T^0 (\xi )d\xi }  =  - \int\limits_{h_0^0}^{h_u^0} {G_{3\xi } (x,t;\xi ,0)T^0 (\xi )d\xi }  \hfill \\
   &=& -G_3 (x,t;h_u ^0,0)T^0 (h_u^0-)+G_3 (x,t;h_0 ^0,0)T^0 (h_0^0-) + \int\limits_{h_0^0}^{h_u^0} {G_3 (x,t;\xi ,0)T_\xi ^0 (\xi )d\xi }.
\eqq
Using the compatible conditions $T^0(h_u^0-)=0$ and $T_0(0)=T^0(h_0^0)$ and replacing $T_0'(\tau)$ by $v_0(\tau)$, we find from (\ref{ieTFb}) that
\bq
  T_x (x,t) = D_I \int\limits_0^t {G_{3x} (x,t;h_u (\tau ),\tau )v_3 (\tau )d\tau }  - D_I \int\limits_0^t {G_{3x} (x,t;h_0 (\tau ),\tau )v_2 (\tau )d\tau }\nn  \hfill \\
 \h\h  - \int\limits_0^t {G_3 (x,t;h_0 (\tau ),\tau )v_0(\tau )d\tau }  + \int\limits_{h_0 ^0}^{h_u^0} {G_3 (x,t;\xi ,0)T_\xi ^0 (\xi )d\xi }\label{TxFb}
\eq
for all $h_u(t)>x>h_0(t)$. Taking $x\to h_0(t)+$ and $x\to h_u(t)-$ in (\ref{TxFb}) and using Lemma \ref{jumpLemma} again we obtain the last two integral equations
\bq
  v_2 (t) &=& 2D_I \int\limits_0^t {G_{3x} (h_0 (t),t;h_u (\tau ),\tau )v_3 (\tau )d\tau }  - 2D_I \int\limits_0^t {G_{3x} (h_0 (t),t;h_0 (\tau ),\tau )v_2 (\tau )d\tau } \nn\hfill \\
&~& - 2\int\limits_0^t {G_3 (h_0 (t),t;h_0 (\tau ),\tau )v_0 (\tau )d\tau }  + 2\int\limits_{h_0^0}^{h_u^0} {G_3 (h_0 (t),t;\xi ,0)T_\xi ^0 (\xi )d\xi } ,  \label{iev2}
\eq
and
\bq
  v_3 (t) &=& 2D_I \int\limits_0^t {G_{3x} (h_u (t),t;h_u (\tau ),\tau )v_3 (\tau )d\tau }  - 2D_I \int\limits_0^t {G_{3x} (h_u (t),t;h_0 (\tau ),\tau )v_2 (\tau )d\tau }\nn  \hfill \\
&~&  - 2\int\limits_0^t {G_3 (h_u (t),t;h_0 (\tau ),\tau )v_0(\tau )d\tau }  + 2\int\limits_{h_0^0}^{h_u^0} {G_3 (h_u (t),t;\xi ,0)T_\xi ^0 (\xi )d\xi } . \label{iev3}
\eq

Note that due to  the  interface condition (\ref{Fb-Oc1}) and (\ref{Fb-Wa}), one has
\bq
h_0 (t) = h_0^0 + \widetilde\lambda _I \int\limits_0^t {v_2 (\tau )d\tau }  - \widetilde\lambda _O \int\limits_0^t {v_1 (\tau )d\tau } .\label{h0v1v2}
\eq
and
\bq
h_u (t) = h_u^0 + \widetilde\lambda _I \int\limits_0^t {v_3(\tau )d\tau }.\label{huv3}
\eq
Therefore the equations (\ref{iev0}), (\ref{iev1}), (\ref{iev2}) and (\ref{iev3}) form a system of nonlinear integral equations of the form $v=Pv$, where $v=(v_0,v_1,v_2,v_3)$ and $Pv=(P_0v,P_1v,P_2v,P_3v)$. Note that the time-depending function $v=(v_0,v_1,v_2,v_3)$ is first only continuous in $(0,\sigma)$ by its definition in (\ref{v}). However, because it is a solution of the system of integral equations (\ref{iev0}), (\ref{iev1}), (\ref{iev2}) and (\ref{iev3}), it is indeed continuous in $[0,\sigma]$. We thus have proved the direct part of the following statement.

\begin{proposition}\label{equivalent} The problem of finding a solution $(h_0,h_u,T,S)$ of (\ref{Fb-Oc1})-(\ref{OcS}) in $(0,\sigma)$ is equivalent to the problem of finding a continuous solution $v=(v_0,v_1,v_2,v_3)$ in $[0,\sigma]$ for the system (\ref{iev0}), (\ref{iev1}), (\ref{iev2}) and (\ref{iev3}) such that $h_u(t)>h_0(t)$ in $(0,\sigma)$, where $h_0$ and $h_u$ are defined by (\ref{h0v1v2}) and (\ref{huv3}).
\end{proposition}

We of course just need to prove the converse part of the statement.

\begin{proof} Suppose that $v=(v_0,v_1,v_2,v_3)$ is a continuous solution in $[0,\sigma]$ for the system (\ref{iev0}), (\ref{iev1}), (\ref{iev2}) and (\ref{iev3}) such that $h_u(t)>h_0(t)$ in $(0,\sigma)$, $h_0$ and $h_u$ are defined by (\ref{h0v1v2}) and (\ref{huv3}). We shall now recover the solution $(h_0,h_u,T,S)$ of (\ref{Fb-Oc1})-(\ref{OcS}).

Since $h_0$ and $h_u$ have already been defined by (\ref{h0v1v2}) and (\ref{huv3}), it remains to determine $T$ and $S$. We define $T_0(t)$ by
\bq
T_0 (t) = T^0 (h_0^0) + \int\limits_0^t {v_0(\tau )d\tau }. \label{T0}
\eq
and determine naturally $T$ by (\ref{ieTOc}) and (\ref{ieTFb}). Similarly, we define $S_0$ by (\ref{v0S0}) and determine $S$ by (\ref{ieS}). We need to prove that $(h_0,h_u,T,S)$ satisfies (C1)-(C4).

{\bf Step 1.} We first prove (C1)-(C3) except the behavior of $T$ and $S$ at the interfaces $x=h_0(t)$ and $x=h_u(t)-$. In fact, the condition (C1) follows the definition of $h_0$ in (\ref{h0v1v2}), the definition of $h_u$ in (\ref{huv3}), and the condition $h_u(t)>h_0(t)$ in definition of the solution $v=(v_0,v_1,v_2,v_3)$.

Due to definition of $S$ in the integral forms (\ref{ieS}), $S_t$ and $S_{xx}$ are continuous in $t\in (0,\sigma)$, $x\in (-\infty,h_0(t))$. Moreover, getting $t\to 0+$ in (\ref{ieS}) and using $\mathop {\lim }\limits_{t \to 0 + } G_1(x,t;\xi ,0) = \delta (x - \xi )$, we obtain the initial condition $S(x,0+) = S^0 (x)$ for all $x\in (- \infty,h_0^0)$.

Similarly, due to definition of $T$ in (\ref{ieTOc}) and (\ref{ieTFb}), $T_t$ and $T_{xx}$ are continuous in $t\in (0,\sigma)$, $x\in (-\infty,h_0(t)) \cup (h_0(t),h_u(t))$, and $T(x,0+) = T^0 (x)$ for all $x\in (-\infty,h_0^0) \cup (h_0^0,h_u^0)$.

{\bf Step 2.} We next check three diffusion equations (\ref{Fb})-(\ref{OcS}). To prove (\ref{OcS}) from the definition of $S$ in (\ref{ieS}), we simply verify that that each term in the right-hand side of (\ref{ieS}) is a homogeneous solution of the operator $(\partial /\partial t - D\partial ^2 /\partial x^2)$ in $t\in (0,\sigma)$, $x\in (-\infty,h_0(t))$. In fact, for the second term one it follows from the properties of Green function that
\bqq
&~& \left( {\frac{\partial }
{{\partial t}} - D\frac{{\partial ^2 }}
{{\partial x^2 }}} \right)\left( {\int\limits_{-\infty}^{h_0^0} {G_1 (x,t;\xi ,0)S^0 (\xi )d\xi } } \right) \hfill \\
 &=& \int\limits_{-\infty}^{h_0^0}  {\left( {\frac{\partial }
{{\partial t}} - D\frac{{\partial ^2 }}
{{\partial x^2 }}} \right)\left( {G_1 (x,t;\xi ,0)} \right)S^0 (\xi )d\xi }  = 0.
\eqq
For the first term, we have to be more careful because in general differentiating with respect to $t$ a function of the form $t \mapsto \int\limits_0^t {K(t,\tau )d\tau }$ may cause a jump,
\[
\frac{d}
{{dt}}\left( {\int\limits_0^t {K(t,\tau )d\tau } } \right) = \mathop {\lim }\limits_{\tau  \to t - } K(t,\tau ) + \int\limits_0^t {K_t (t,\tau )d\tau } .
\]
However, in this case the jump $\mathop {\lim }\limits_{\tau  \to t - } K(t,\tau )$ vanishes since
\[
\mathop {\lim }\limits_{\tau  \to t - } G_1 (x,t;\xi ,\tau ) = \mathop {\lim }\limits_{\tau  \to t - } G_{1\xi } (x,t;\xi ,\tau ) = 0,~~\forall x \ne \xi.
\]
Therefore, the first term can be treat similarly to the second term. Thus (\ref{OcS}) holds. The proofs for (\ref{Fb}) and (\ref{OcT}) are the same.

{\bf Step 3.} We now prove that $S(h_0(t)-,t)=S_0(t)$, which in
particular accomplish (C3), and verify the Stefan condition
(\ref{Fb-Oc2}).

In fact, (\ref{ieS0}) holds due to (\ref{v0S0}) and (\ref{iev0}). On the other hand, taking $x\to h_0(t)-$ in (\ref{ieS}) and using Lemma \ref{jumpLemma} one has
\bqq
  S(h_0 (t) - ,t) &=& \frac{1}
{2}S_0 (t) - D\int\limits_0^t {G_{1\xi } (h_0 (t),t;h_0 (\tau ),\tau )S_0 (\tau )d\tau }  \hfill \\
&~& + \int\limits_{ - \infty }^{h_0^0} {G_1 (h_0 (t),t;\xi ,0)S^0(\xi)d\xi } .
\eqq
It follows from the latter equation and (\ref{ieS0}) that $S(h_0 (t)- ,t)=S_0(t)$.

Applying Lemma \ref{mainLemma} to equation (\ref{OcS}) for $h_0(t)>x>-\infty$ and using $S(h_0(t)-,t)=S_0(t)$, we have
\bqq
  S(x,t) &=& \int\limits_0^t {G_1 (x,t;h_0 (\tau ),\tau )\left[ {DS_x (h_0 (\tau ) - ,\tau ) + S_0 (\tau )h_0 '(\tau )} \right]d\tau }  \hfill \\
&~&- D\int\limits_0^t {G_{1\xi } (x,t;h_0 (\tau ),\tau )S_0 (\tau )d\tau }  + \int\limits_{ - \infty }^{h_0^0}  {G_1 (x,t;\xi ,0)S^0(\xi)d\xi } .
\eqq
Comparing the latter equation to the original definition of $S$ in (\ref{ieS}), we obtain
\bq
\int\limits_0^t {G_1 (x,t;h_0 (\tau ),\tau )\Psi _1 (\tau )d\tau }  = 0,\h x\in (-\infty,h_0(t)),\label{AF}
\eq
where $\Psi _1 (t) = DS_x (h_0 (t) - ,t) + S_0 (t)h_0 '(t)$. We shall deduce from (\ref{AF}) that $\Psi _1 \equiv 0$, which is equivalent to (\ref{Fb-Oc2}). Indeed, differentiating (\ref{AF}) with respect to $x$, then taking $x\to h_0(t)-$ and using Lemma \ref{jumpLemma} we get
\bqq
\Psi _1 (t) =  - 2\int\limits_0^t {G_{1x} (h_0 (t),t;h_0 (\tau ),\tau )\Psi _1 (\tau )d\tau }.
\eqq
Note that this is a linear Volterra integral equation of the second kind and
\[
\left| {G_{1x} (h_0 (t),t;h_0 (\tau ),\tau )} \right| \leqslant \frac{{const.}}
{{\sqrt {t - \tau } }}.
\]
Therefore, it follows from Lemma \ref{linearVolterra} that $\Psi _1 (t)=0$.

{\bf Step 4.} We prove that $T_x(h_0(t)-,t)=v_1(t)$, $T_x(h_0(t)+,t)=v_2(t)$ and $T_x(h_u(t)-,t)=v_3(t)$, which in particular imply the Stefan conditions (\ref{Fb-Oc1}) and (\ref{Fb-Wa}). Note that we have already had from (\ref{T0}) that $v_0(t)=T_0'(t)$ and $T_0(0)=T^0(h_0^0)$.

To prove $T_x(h_u(t)-,t)=v_3(t)$ from the definition of $T$ in
(\ref{ieTFb}), we use the same process of getting (\ref{iev3}) from
(\ref{ieTFb}). In fact, differentiating both sides of (\ref{ieTFb})
with respect to $x$ and then taking $x\to h_u(t)-$ we have \bqq
  T_x (h_u (t) - ,t) &=& \frac{1}
{2}v_3 (t) + D_I \int\limits_0^t {G_{3x} (h_u (t),t;h_u (\tau ),\tau )v_3 (\tau )d\tau }  \hfill \\
&~& - D_I \int\limits_0^t {G_{3x} (h_u (t),t;h_0 (\tau ),\tau )v_2 (\tau )d\tau }  \hfill \\
&~&  - \int\limits_0^t {G_3 (h_u (t),t;h_0 (\tau ),\tau )T_0 '(\tau )d\tau }  + \int\limits_{h_0^0}^{h_u^0} {G_3 (h_u (t),t;\xi ,0)T_\xi ^0 (\xi )d\xi }.
\eqq
Comparing the latter equation to (\ref{iev3}) and using $v_0(\tau)=T_0'(\tau)$, we obtain $T_x (h_u (t) - ,t)=v_3(t)$. By the same way we find that $T_x (h_0 (t) + ,t)=v_2(t)$ and  $T_x (h_0 (t) - ,t)=v_1(t)$.

{\bf Step 5.} Finally, we verify that $T(h_0(t),t)=T_0(t)$ and $T(h_u(t)-,t)=0$ to accomplish (C2) and (C4).

Applying Lemma \ref{mainLemma} to equation (\ref{OcT}) for $h_0(t)>x>-\infty$ one has
\bqq
  T(x,t) &=& \int\limits_0^t {G_2 (x,t;h_0 (\tau ),\tau )\left[ {D_O T_x (h_0 (\tau ) - ,\tau ) + T(h_0 (\tau ) - ,\tau )h_0 '(\tau )} \right]d\tau }  \hfill \\
&~& - D_O \int\limits_0^t {G_{2\xi } (x,t;h_0 (\tau ),\tau )T(h_0 (\tau ) - ,\tau )d\tau }  + \int\limits_{ - \infty }^{h_0^0} {G_2 (x,t;\xi ,0)T^0 (\xi )d\xi }.
\eqq
Comparing the latter equation to the original definition of $T$ in (\ref{ieTOc}) and using $T_x (h_0 (t) - ,t)=v_1(t)$, we obtain
\bq
\int\limits_0^t {\left[ {G_2 (x,t;h_0 (\tau ),\tau )h_0 '(\tau ) - D_O G_{2\xi } (x,t;h_0 (\tau ),\tau )} \right]\Psi_2 (\tau )d\tau }  = 0
,\label{Dir1}
\eq
where $\Psi_2 (t) = T(h_0 (t) - ,t) - T_0 (t)$. We now prove $\Psi_2\equiv 0$ by using the same technique of dealing $\Psi_1\equiv 0$. Taking $x\to h_u(t)-$ in (\ref{Dir1}) and using Lemma \ref{jumpLemma}, we find that
$$
\Psi_2 (t) = 2\int\limits_0^t {\left[ {G_2 (x,t;h_0 (\tau ),\tau )h_0 '(\tau ) - D_O G_{2\xi } (x,t;h_0 (\tau ),\tau )} \right]\Psi_2 (\tau )d\tau }. \label{Dir1a}
$$
Note that this is a linear Volterra integral equation of the second
kind and
$$
\left| {G_2 (x,t;h_0 (\tau ),\tau )h_0 '(\tau ) - D_O G_{2\xi } (x,t;h_0 (\tau ),\tau )} \right| \leqslant \frac{{const.}}
{{\sqrt {t - \tau } }}.
$$
Therefore, it follows from Lemma \ref{linearVolterra} that $\Psi_2 \equiv 0$. Thus $T(h_0 (t) - ,t)=T_0 (t)$.

We shall use the same argument, in fact a little more complicated one, to deduce that $T(h_0 (t)+ ,t)=T_0 (t)$ and $T(h_u (t),t)=0$. Applying Lemma \ref{mainLemma} to equation (\ref{Fb}) for $h_u(t)>x>h_0(t)$ and comparing to the original definition of $T$ in (\ref{ieTFb}), one has
\bq
&~& \int\limits_0^t {\left[ {G_3 (x,t;h_0 (\tau ),\tau )h_0 '(\tau ) - D_I G_{3\xi } (x,t;h_0 (\tau ),\tau )} \right]\Psi _3 (\tau )d\tau } \nn \hfill \\
   &=& \int\limits_0^t {\left[ {G_3 (x,t;h_u (\tau ),\tau )h_u '(\tau ) - D_I G_{3\xi } (x,t;h_u (\tau ),\tau )} \right]\Psi _4 (\tau )d\tau } ,\label{Dir23}
\eq
where $\Psi _3 (\tau ) = T(h_0 (t) + ,t) - T_0 (t)$ and $\Psi_4 (t) = T(h_u (t)-,t)$. Taking $x\to h_0(t)+$ and $x\to h_u(t)-$ in (\ref{Dir23}), we have
\bq
  \Psi _3 (t) &=& 2\int\limits_0^t {\left[ {G_3 (h_0 (t),t;h_0 (\tau ),\tau )h_0 '(\tau ) - D_I G_{3\xi } (h_0 (t),t;h_0 (\tau ),\tau )} \right]\Psi _3 (\tau )d\tau } \nn \hfill \\
&-&2\int\limits_0^t {\left[ {G_3 (h_0 (t),t;h_u (\tau ),\tau )h_u '(\tau ) - D_I G_{3\xi } (h_0 (t),t;h_u (\tau ),\tau )} \right]\Psi _4 (\tau )d\tau }, \label{Dir2}
\eq
and
\bq
  -\Psi _4 (t) &=&  -2 \int\limits_0^t {\left[ {G_3 (h_u (t),t;h_0 (\tau ),\tau )h_0 '(\tau ) - D_I G_{3\xi } (h_u (t),t;h_0 (\tau ),\tau )} \right]\Psi _3 (\tau )d\tau }\nn  \hfill \\
&+&2 \int\limits_0^t {\left[ {G_3 (h_u (t),t;h_u (\tau ),\tau )h_u '(\tau ) - D_I G_{3\xi } (h_u (t),t;h_u (\tau ),\tau )} \right]\Psi _4 (\tau )d\tau } .  \label{Dir3}
\eq
Note that equations $(\ref{Dir2})$ and $(\ref{Dir3})$ form a system of linear Volterra integral equations of the second kind and
\bqq
\left| {G_3 (h_u (t),t;h_0 (\tau ),\tau )h_0 '(\tau ) - D_I G_{3\xi } (h_u (t),t;h_0 (\tau ),\tau )} \right| &\leqslant& \frac{{const.}}
{{\sqrt {t - \tau } }},\hfill\\
\left| {G_3 (h_0 (t),t;h_u (\tau ),\tau )h_u '(\tau ) - D_I G_{3\xi } (h_0 (t),t;h_u (\tau ),\tau )} \right| &\leqslant& \frac{{const.}}
{{\sqrt {t - \tau } }}.
\eqq
We therefore deduce from Lemma \ref{linearVolterra} that $\Psi_3\equiv  \Psi_4 \equiv 0$. Thus $T(h_u(t)-,t)=0$, i.e. (\ref{Fb-Wa1}), and $T(h_0(t),t)=T_0(t)$, which accomplishs (C2) and also implies (\ref{T0S0}) in viewing (\ref{v0S0}). The proof is completed.
\end{proof}

\begin{remark} We learned the technique of proving $\Psi_1\equiv 0$ by (\ref{AF}) from Friedman \cite{F}. However, there is a minor mistake in \cite{F}. The formula (1.28) in page 221, i.e.
$$
\int\limits_0^t {G_\xi  (x,t;s(\tau ),\tau )u(s(\tau ),\tau )d\tau }  = 0,
$$
should be
$$
\int\limits_0^t {\left[ {G_\xi  (x,t;s(\tau ),\tau ) - G(x,t;s(\tau ),\tau )s'(\tau )} \right]u(s(\tau ),\tau )d\tau }  = 0.
$$
Fortunately, the conclusion $u(s(t ),t )=0$ still holds by the same argument.
\end{remark}

\section{Solving integral equations}\label{uniqueness4}

\h Now we accomplish the proof of Theorem \ref{existenceanduniqueness}. We want to prove the system (\ref{iev0}), (\ref{iev1}), (\ref{iev2}) and (\ref{iev3}) has a unique local solution $v=(v_0,v_1,v_2,v_3)$, and this solution can be extended uniquely by prolongation whenever the condition $h_u(t)>h_0(t)$ still holds. Note that this system is of the form $v=Pv:=\left( {P_0 v, P_1v,P_2 v, P_3 v}\right)$, where $P_0v,P_1v, P_2v$ and $P_3 v$ are the right-hand sides of (\ref{iev0}), (\ref{iev1}), (\ref{iev2}) and (\ref{iev3}) respectively.

Denote by $C(\sigma,M)$ the space of functions $v=(v_0,v_1,v_2,v_3)$ continuous in $[0,\sigma]$ and
\[
\left\| v \right\|_{[0,\sigma ]} : = \mathop {\max }\limits_{t \in \left[ {0,\sigma } \right],j = \overline {0,3}}|v_j(t)| \leqslant M.
\]
When $v=0$ then $P(0)$ depends only on the initial conditions $(h_0^0,h_u^0,T^0,S^0)$ and is continuous on $[0,\infty)$. Fix $M>0$ large enough, says $M \geqslant 2\left\| P(0)\right\|_{[0,1]}$.

\begin{proposition}\label{Pcontraction} There exists $\sigma>0$ depending on $(h_0^0,h_u^0,T^0,S^0)$ such that $P$ maps $C(t_0,M)$ into itself and is a contraction.
\end{proposition}

Note that $C(\sigma,M)$ is a complete metric space. Therefore it follows from Proposition \ref{Pcontraction} and the Contraction Mapping Principle that $P$ has a unique fixed point $v$ in $C(\sigma,M)$. It gives also the unique solution $(h_0,h_u,T,S)$ of (\ref{Fb-Oc1})-(\ref{OcS}) in $(0,\sigma)$ due to Proposition \ref{equivalent}.

\begin{proof} First at all, let us consider $\sigma>0$ small enough, says
\[
\sigma  \leqslant \frac{{h_u^0  - h_0^0 }}
{{2(2\widetilde\lambda _I  + \widetilde\lambda _O )M}}.
\]
Then for $v\in C(\sigma,M)$, it follows from (\ref{h0v1v2})-(\ref{huv3}) that
\bq
  h_u (t) - h_0 (\tau) &\ge& h_u^0-h_0^0-\left| {h_u (t) - h_u (0)} \right|-\left| {h_0 (\tau) - h_0 (0)} \right|
\nn\hfill\\
  &\ge& h_u^0-h_0^0-\widetilde\lambda _IMt-(\widetilde\lambda _O+\widetilde\lambda _I)M\tau\ge \frac{{h_u^0- h_0^0}}
{2} > 0\label{huh0}
\eq
for all $0\le 0\le t\le \sigma$. In particular $h_u(t)>h_0(t)$ for all $t\in [0,\sigma]$ and hence $P$ is well-defined on $C(\sigma,M)$. We now estimate $\left\| {Pv - P\widetilde v} \right\|_{[0,\sigma]}$ for $v,\widetilde v\in C(\sigma,M)$. In what follows, denote by $\widetilde h_u$ and $\widetilde h_0$ the functions given by (\ref{h0v1v2})-(\ref{huv3}) where $v$ is replaced by $\widetilde v$. For simplicity, we shall always denote by $C_0>0$ an universal constant depending only on the initial conditions $(h_0^0,h_u^0,T^0,S^0)$. We shall go to the details for $P_3$, and $P_0, P_1,P_2$ can be treated by the same way. Recall that
\bqq
  (P_3v)(t) &=& 2D_I \int\limits_0^t {G_{3x} (h_u (t),t;h_u (\tau ),\tau )v_3 (\tau )d\tau }  - 2D_I \int\limits_0^t {G_{3x} (h_u (t),t;h_0 (\tau ),\tau )v_2 (\tau )d\tau }\nn  \hfill \\
&~&  - 2\int\limits_0^t {G_3 (h_u (t),t;h_0 (\tau ),\tau )v_0(\tau )d\tau }  + 2\int\limits_{h_0^0}^{h_u^0} {G_3 (h_u (t),t;\xi ,0)T_\xi ^0 (\xi )d\xi } .
\eqq
We need some preliminary estimates related to the Green's function $G_3$.

\begin{lemma}\label{prepareG} For all $v,\widetilde v\in C(\sigma,M)$, $0\le \tau<t\le \sigma$ and $\xi\in \mathbb{R}$, we have
\bq
  &~&\left| {G_{3x} (h_u (t),t;h_u (\tau ),\tau ) - G_{3x} (\widetilde h_u (t),t;\widetilde h_u (\tau ),\tau )} \right| \leqslant \frac{{C_0}}
{{\sqrt {t - \tau } }}\left\| {v - \widetilde v} \right\|_{[0,\sigma]} , \label{G31}\hfill \\
  &~&\left| {G_3 (h_u (t),t;h_0 (\tau ),\tau ) - G_3 (\widetilde h_u (t),t;\widetilde h_0 (\tau ),\tau )} \right| \leqslant \frac{{C_0}}
{{\sqrt {t - \tau } }}\left\| {v - \widetilde v} \right\|_{[0,\sigma]} , \label{G32} \hfill \\
 &~& \left| {G_{3x} (h_u (t),t;h_0 (\tau ),\tau ) - G_{3x} (\widetilde h_u (t),t;\widetilde h_0 (\tau ),\tau )} \right| \leqslant \frac{{C_0}}
{{\sqrt {t - \tau } }}\left\| {v - \widetilde v} \right\|_{[0,\sigma]} , \label{G33} \hfill \\
 &~& \left| {G_{3x} (h_u (t),t;\xi ,\tau ) - G_{3x} (\widetilde h_u (t),t;\xi ,\tau )} \right| \leqslant\nn  \hfill \\
   &~&\h\leqslant C_0 \left[ {\exp \left( { - \frac{{(h_u (t) - \xi )^2 }}
{{8D_I t}}} \right) + \exp \left( { - \frac{{(h_u (t) - \xi )^2 }}
{{8D_I t}}} \right)} \right]\left\| {v - \widetilde v} \right\|_{[0,\sigma ]}\label{G34} .
\eq
\end{lemma}

\begin{proof} Considering
\[
G_{3x} (h_u (t),t;h_u (\tau ),\tau ) = -\frac{{h_u (t) - h_u (\tau )}}
{{4\sqrt {\pi D_I^3 (t - \tau )^3 } }}\exp \left( { - \frac{{(h_u (t) - h_u (\tau ))^2 }}
{{4D_I (t - \tau )}}} \right)
\]
as a function of one-variable $(h_u (t) - h_u (\tau ))$ and using Lagrange formula one has
\bqq
&~&  G_{3x} (h_u (t),t;h_u (\tau ),\tau ) - G_{3x} (\widetilde h_{u} (t),t;\widetilde h_{u} (\tau ),\tau ) \hfill \\
   &=& -\frac{{[ {\left( {h_u (t) - h_u (\tau )} \right) - ( {\widetilde h_{u} (t) - \widetilde h_{u} (\tau )} )} ]}}
{{4\sqrt {\pi D_I^3 (t - \tau )^3 } }}\exp \left( { - \frac{{\theta ^2 }}
{{4D_I (t - \tau )}}} \right)\left( {1 - \frac{{2\theta ^2 }}
{{4D_I (t - \tau )}}} \right)
\eqq
for some $\theta$ between $(h_u (t) - h_u (\tau ))$ and $(\widetilde h_u (t) -\widetilde h_u (\tau ))$. Thus (\ref{G31}) follows
\[
|(h_u (t) - h_u (\tau ))-(\widetilde h_u (t) -\widetilde h_u (\tau ))| \leqslant  \widetilde\lambda_I(t - \tau )\left\| {v - \widetilde v} \right\|_{[0,\sigma ]},
\]
and the elementary inequality
\[
e^{ - z} |1 - 2z| \leqslant \frac{{|1 - 2z|}}
{{1 + z}} \leqslant 2,\h\text{where}~ z=\frac{{\theta ^2 }}{{4D_I (t - \tau )}} \geqslant 0.
\]
To prove (\ref{G32}) we shall use the identity
\bqq
&~&G_3 (h_u (t),t;h_0 (\tau ),\tau ) - G_3 (\widetilde h_{u} (t),t;\widetilde h_0 (\tau ),\tau )\hfill\\
&=& -\frac{{[ {(h_u (t) - h_0 (\tau ))-(\widetilde h_{u} (t) - \widetilde h_0(\tau ))}]\theta_1 }}
{{4\sqrt {\pi D_I^3 (t - \tau )^3 } }}\exp \left( { - \frac{{\theta_1^2 }}
{{4D_I (t - \tau )}}} \right)
\eqq
for some $\theta_1$ between $(h_u (t) - h_0 (\tau ))$ and $(\widetilde h_{u} (t) - \widetilde h_0(\tau ))$. Note from (\ref{huh0}) that $|\theta_1|\ge \frac{h_u^0-h_0^0}{2}$. Therefore (\ref{G32}) follows
\[
|(h_u (t) - h_0 (\tau )) - (\widetilde h_u (t) - \widetilde h_0 (\tau ))| \leqslant (2\widetilde\lambda _I  + \widetilde\lambda _O )\sigma \left\| {v - \widetilde v} \right\|_{[0,\sigma ]}
\]
and the elementary inequality
\[
e^{ - z}z  \leqslant 1,\h\text{where}~z = \frac{{\theta _1^2 }}
{{4D_I (t - \tau )}}\ge 0.
\]
We can also prove (\ref{G33}) by the same argument. In fact, we write
\bqq
&~&  G_{3x} (h_u (t),t;h_0 (\tau ),\tau ) - G_{3x} (\widetilde h_{u} (t),t;\widetilde h_0 (\tau ),\tau ) \hfill \\
&=& -\frac{{[(h_u (t) - h_0 (\tau ))-(\widetilde h_{u} (t) - \widetilde h_0(\tau ))]}}
{{4\sqrt {\pi D_I^3 (t - \tau )^3 } }}\exp \left( { - \frac{{\theta_2^2 }}
{{4D_I (t - \tau )}}} \right)\left( {1 - \frac{{2\theta_2^2 }}
{{4D_I (t - \tau )}}} \right)
\eqq
for some $\theta_2$ between $(h_u (t) - h_0 (\tau ))$ and $(\widetilde h_{u} (t) - \widetilde h_0(\tau ))$ and use
\[
e^{ - z} z|1 - 2z|\leqslant \frac{{z|1 - 2z|}}
{{1+z + \frac{1}
{2}z^2 }} \leqslant 4,\h\text{where}~z=\frac{{\theta _2^2 }}
{{4D_I (t - \tau )}} \geqslant 0.
\]
Finally, for the last inequality (\ref{G34}) we write
$$
G_3 (h_u (t),t;\xi,0 ) - G_3 (\widetilde h_{u} (t),t;\xi,0)= -\frac{{[ {h_u (t) -\widetilde h_{u} (t)}]\theta_3 }}
{{4\sqrt {\pi D_I^3 t^3 } }}\exp \left( { - \frac{{\theta_3^2 }}
{{4D_I t}}} \right)
$$
for some $\theta_3$ between $(h_u(t) -\xi)$ and $(\widetilde h_{u} (t) -\xi)$. Note that
$$
\exp \left( { - \frac{{\theta _3^2 }}
{{8D_I t}}} \right) \leqslant \max \left\{ {\exp \left( { - \frac{{(h_u (t) - \xi )^2 }}
{{8D_I t}}} \right),\exp \left( { - \frac{{(\widetilde h_u (t) - \xi )^2 }}
{{8D_I t}}} \right)} \right\} .
$$
Thus (\ref{G34}) follows
$$
  e^{- z} \sqrt z  \leqslant \frac{{\sqrt z }}
{{1 + z}} \leqslant \frac{1}
{2},\h\text{where}~z = \frac{{\theta _3^2 }}
{{8D_I t}} \geqslant 0,
$$
and the estimate $| {h_u (t) - \widetilde h_u (t)} | \leqslant \widetilde\lambda _I t\left\| {v - \widetilde v} \right\|_{[0,\sigma]}$.
\end{proof}

We now already to estimate $||P_3v-P_3\widetilde v||_{[0,\sigma]}$ for $v,\widetilde v\in C(\sigma,M)$. Recall that $P_3$ is the sum of four terms
\bqq
  (P_3v)(t) &=& 2D_I \int\limits_0^t {G_{3x} (h_u (t),t;h_u (\tau ),\tau )v_3 (\tau )d\tau }  - 2D_I \int\limits_0^t {G_{3x} (h_u (t),t;h_0 (\tau ),\tau )v_2 (\tau )d\tau }\nn  \hfill \\
&~&  - 2\int\limits_0^t {G_3 (h_u (t),t;h_0 (\tau ),\tau )v_0(\tau )d\tau }  + 2\int\limits_{h_0^0}^{h_u^0} {G_3 (h_u (t),t;\xi ,0)T_\xi ^0 (\xi )d\xi } .
\eqq
For the first term of $P_3$, using (\ref{G31}) we have
\bqq
&~&| {G_{3x} (h_u (t),t;h_u (\tau ),\tau )v_3 (\tau ) - G_{3x} (\widetilde h_u (t),t;\widetilde h_u (\tau ),\tau )\widetilde v_3 (\tau )} | \hfill \\
 &\leqslant& | {G_{3x} (h_u (t),t;h_u (\tau ),\tau ) - G_{3x} (\widetilde h_u (t),t;\widetilde h_u (\tau ),\tau )} |.| {v_3 (\tau )} | \hfill \\
&~&~  + | {G_{3x} (\widetilde h_u (t),t;\widetilde h_u (\tau ),\tau )}|.| {v_3 (\tau ) - \widetilde v_3 (\tau )}| \hfill \\
 &\le&\frac{{C_0 M}}
{{\sqrt {t - \tau } }}\left\| {v - \widetilde v} \right\|_{[0,\sigma ]} + \frac{{C_0}}
{{\sqrt {t - \tau } }} \left\| {v - \widetilde v} \right\|_{[0,\sigma ]}=\frac{{C_0(M+1)}}
{{\sqrt {t - \tau } }} \left\| {v - \widetilde v} \right\|_{[0,\sigma ]}
\eqq
for all $0\le \tau\le t\le \sigma$. Therefore
\bqq
&~&\int\limits_0^t {\left| {G_{3x} (h_u (t),t;h_u (\tau ),\tau )v_3 (\tau ) - G_{3x} (\widetilde h_u (t),t;\widetilde h_u (\tau ),\tau )\widetilde v_3 (\tau )} \right|d\tau }\nn  \hfill \\
&\le& C_0 (M+1)\sqrt {\sigma} \left\| {v - \widetilde v} \right\|_{[0,\sigma]},\h t\in [0,\sigma].  \label{EP1}
\eqq
We have also the same estimate for the second term and the third term of $P_3$ due to (\ref{G32}) and (\ref{G33}). For the last term of $P_3$, it follows from (\ref{G34}) and the boundness of $T^0_x$ that
\bqq
&~&\int\limits_{h_0^0}^{h_u^0} {\left| {G_3 (h_u (t),t;\xi ,0) - G_3 (\widetilde h_u (t),t;\xi ,0)} \right|T_\xi ^0 (\xi )d\xi }  \hfill \\
 &\leqslant& C_0 \left\| {v - \widetilde v} \right\|_{[0,\sigma ]} \int\limits_{h_0^0}^{h_u^0} {\left[ {\exp \left( { - \frac{{( h_u (t)-\xi )^2 }}
{{8D_I t}}} \right) + \exp \left( { - \frac{{(\widetilde h_u (t)-\xi )^2 }}
{{8D_I t}}} \right)} \right]d\xi } .
\eqq
Note that by changing of variables
\[
\begin{gathered}
 ~~\int\limits_{ - \infty }^\infty  {\left[ {\exp \left( { - \frac{{(h_u (t)-\xi )^2 }}
{{8D_It}}} \right) + \exp \left( { - \frac{{(\widetilde h_{u} (t)-\xi )^2 }}
{{8D_I t}}} \right)} \right]d\xi }  \hfill \\
   = 2\int\limits_{ - \infty }^\infty  {\exp \left( { - \frac{{\xi ^2 }}
{{8D_I t}}} \right)d\xi }  = 2\sqrt {\pi 8D_It} . \hfill \\
\end{gathered}
\]
Thus
$$\int\limits_{h_0 (0)}^{h_u (0)} {\left| {G_3 (h_u (t),t;\xi ,0) - G_3 (\widetilde h_u (t),t;\xi ,0)} \right|T_\xi ^0 (\xi )d\xi }\leqslant 2C_0\sqrt {\pi 8D_I \sigma } \left\| {v - \widetilde v} \right\|_{[0,\sigma]} .
$$
In summary, we have
\[
\left\| {P_3 v - P_3 \widetilde v} \right\|_{[0,\sigma]}  \leqslant C_0 \sqrt \sigma \left\| {v - \widetilde v} \right\|_{[0,\sigma]} ,
\]
where $C_0$ stands for a constant depending only on the initial conditions $(h_0^0,h_u^0,T^0,S^0)$. We have also the same estimates for $P_0,P_1,P_2$. Thus for $\sigma>0$ small enough, says $\sigma\le \frac{1}{4C_0^2}$, one has
\bq
\left\| {Pv - P\widetilde v} \right\|_{[0,\sigma]}  \leqslant \frac{1}
{2}\left\| {v - \widetilde v} \right\|_{[0,\sigma]},~\forall v,\widetilde v\in C(\sigma,M). \label{contraction}
\eq
Note that $M \geqslant 2\left\| P(0)\right\|_{[0,1]}\ge 2\left\| P(0)\right\|_{[0,\sigma]}$. Thus it follows from (\ref{contraction}) that $P$ maps from $C(\sigma,M)$ into itself and is a contraction.
\end{proof}

To finish the proof of Theorem \ref{existenceanduniqueness}, we prove that the solution can be extended uniquely by prolongation whenever $h_u(\sigma)>h_0(\sigma)$. Proposition \ref{equivalent} allows us consider the problem of finding $v=(v_0,v_1,v_2,v_3)$ instead of the problem of finding $(h_0,h_u,T,S)$. Assume  that $v=(v_0,v_1,v_2,v_3)$ is a solution of the system (\ref{iev0}), (\ref{iev1}), (\ref{iev2}) and (\ref{iev3}) in $[0,\sigma]$ such that $h_u(t)>h_0(t)$ for all $t\in [0,\sigma]$. From such solution $v$, we can extend naturally $(h_0,h_u,T,S)$ in $[0,\sigma]$ by using the scheme of the proof of Proposition \ref{equivalent}. We now want to check that $(h_0(\sigma),h_u(\sigma),T(.,\sigma),S(.,\sigma))$ satisfy (H1)-(H3).

In fact, the condition (H1) is automatically satisfied. For (H2),
the continuity of $T_x(.,\sigma)$ in $(-\infty,h_0(\sigma))\cup
(h_0(\sigma),h_u(\sigma))$ is guaranteed by (\ref{TxOc}) and
(\ref{TxFb}), which still hold for $t=\sigma$. The behavior of $T$
at the interface is also insured by the relationships $T_x (h_0
(\sigma) - ,\sigma)=v_1 (\sigma)$, $T_x (h_0 (\sigma) + ,\sigma) =
v_2 (\sigma)$, $T_x (h_u (\sigma) - ,\sigma) = v_3 (\sigma)$,
$T(h_0(\sigma),\sigma)=T_0(\sigma)$ and $T(h_u(\sigma)-,\sigma)=0$.
Moreover, we deduce from (\ref{TxOc}) that
\[
\mathop {\lim \sup }\limits_{x \to  - \infty } \left| {T_x (x,\sigma)} \right| = \mathop {\lim \sup }\limits_{x \to  - \infty } \left| {\int\limits_{ - \infty }^{h_0 (0)} {G_2 (x,t;\xi ,0)T_\xi ^0 (\xi )d\xi } } \right| \leqslant \mathop {\sup }\limits_{h_0 (0) > \xi  >  - \infty } \left| {T_\xi ^0 (\xi )} \right|
\]
and the boundness of $T_x$ follows. Here we have used a property of
the Green's function
\[
\int\limits_{ - \infty }^\infty  {\left| {G_2 (x,t;\xi ,0)} \right|d\xi }  = 1.
\]
Thus (H2) is indeed true. Similarly, (H3) holds.

By considering $(h_0(\sigma),h_u(\sigma), T(.,\sigma),S(.,\sigma))$ as the new initial condition, we use Proposition \ref{Pcontraction} to extend uniquely the solution $v=(v_0,v_1,v_2,v_3)$ of the system (\ref{iev0}), (\ref{iev1}), (\ref{iev2}) and (\ref{iev3}) in $[0,\sigma^*]$ such that $h_u(t)>h_0(t)$ for all $t\in (0,\sigma^*)$ for some $\sigma^*>\sigma$. Finally, it follows from Proposition \ref{equivalent} that the extended solution $v$ gives the extended solution $(h_0,h_u,T,S)$ of the system (\ref{Fb-Oc1})-(\ref{OcS}) in $(0,\sigma^*)$ corresponding to the initial conditions $(h_0^0, h_u^0,T^0,S^0)$. The proof is completed.

\begin{remark} We may also replace the infinite domain $h_u(0)>x>-\infty$ of the initial conditions by a finite interval $h_u^0>x>L$ where $L$ is a deep point in the ocean satisfying $L<h_0(t)$ for $0\le t \le \sigma$. However, in this case it is necessary to require the history information $(T(L,.),S(L,.))$, $0<t<\sigma$, at the point $L$. We can find the same result by using the Green's function for the half-plane $x>L$,
\[
K_j (x,t;\xi ,\tau ) = G_j (x,t;\xi ,\tau ) - G_j (2L - x,t;\xi ,\tau ),\h j = 1,2,3.
\]
\end{remark}
\text{}\\
{\bf Acknowledgments:} The main part of the paper was accomplished when the second author was doing on the project of his Master thesis in Laboratory MAPMO, University d'Orleans in the summer 2008. He would like to express his heart-felt thanks to everybody of the laboratory for the warm hospitality.

\end{document}